\documentclass[final,1p,times,authoryear,sort]{elsarticle}
\usepackage{amsmath,amsfonts,amssymb,amsthm}
\usepackage[all]{xy}
\usepackage{graphicx,subfig}
\usepackage{aliascnt,hyperref}

\journal{Journal of Symbolic Computation}

\newtheorem{teo}{Theorem}

\def\partet#1#2#3#4{\newaliascnt{#1}{#2}\newtheorem{#1}[#1]{#3}\aliascntresetthe{#1}\providecommand*{#4}{#3}}
\def\parted#1#2#3#4{\newaliascnt{#1}{#2}\newdefinition{#1}[#1]{#3}\aliascntresetthe{#1}\providecommand*{#4}{#3}}
\partet{lema}{teo}{Lemma}{\lemaautorefname}
\partet{prop}{teo}{Proposition}{\propautorefname}
\partet{cor}{teo}{Corollary}{\corautorefname}
\parted{defs}{teo}{Definition}{\defsautorefname}
\parted{ejemplo}{teo}{Example}{\ejemploautorefname}
\parted{obs}{teo}{Remark}{\obsautorefname}
\parted{parr}{teo}{}{\parrautorefname}
\parted{notc}{teo}{Notation}{\notcautorefname}
\parted{ass}{teo}{Assumption}{\assautorefname}
\newproof{dem}{Proof}
\newtheorem*{nteo}{Theorem}

\theoremstyle{definition}

\begin{document}

\begin{frontmatter}

\title{Segmentation of real algebraic plane curves.}

\author[dm]{C\'esar Massri\corref{correspondencia}\fnref{financiado}}
\ead{cmassri@dm.uba.ar}
\address[dm]{Department of Mathematics, FCEN, University of Buenos Aires, Argentina}
\cortext[correspondencia]{Address for correspondence: Department of Mathematics, FCEN,
University of Buenos Aires, Argentina. Postal address: 1428. Phone number: 54-11-4576-3335.}
\fntext[financiado]{The author was fully supported by CONICET, IMAS, Buenos Aires, Argentina}
\author[dc]{Manuel Dubinsky}
\address[dc]{Computer Engineering Department, University of Avellaneda, Argentina}

\begin{abstract}
In this article we give an implementation of
the standard algorithm to segment a real algebraic plane curve defined
implicitly.
Our implementation is efficient and simpler than previous. 
We use global information to count the number
of half-branches at a critical point.
\end{abstract}

\end{frontmatter}

\section*{Introduction.}
\subsection*{Overview of the problem.}
This article is about the computation of a 
geometric-combinatorial description of
a plane curve defined as the zero locus
of a bivariate real polynomial. 

This problem is relevant in computational vision \citep{MR1866897}, 
engineering \citep{MR702260}, architecture \citep{MR1520983} and more.
The algorithms to visualize spacial curves and surfaces rely on
this problem, \citep{MR2532173}.

There are two main types of algorithms 
to produce a geometric-combinatorial description of a plane curve
in polynomial time.
The first one is inspired by 
Morse Theory \citep{MR0163331} and 
uses a sweeping line to detect critical
topological events, such as singularities. 
It goes as follows, see \citep{MR1430842}.
\begin{enumerate}
\item[0.] Ensure that the equation of the curve $g$ is in general position.
\item[1.] Compute the real roots of the discriminant $\Delta_y(g)$, $\{x_i\}$.
\item[2.] Compute the real roots of $g(x_i,-)$, $\{y_{ij}\}$.
\item[3.] For each point $(x_i,y_{ij})$ count the number of half-branches to the left and to the right.
\item[4.] Connect up the points with edges, obeying the branch counts.
\end{enumerate}

The other type of methods relies on subdivision
techniques, \citep{alberti:inria-00175062}. 
The most famous family
of algorithm using this approach is the Marching Cube,
see \citep{Lorensen:1987:MCH:37402.37422}. 

\

We are interested in the first algorithm.
Known implementations of this algorithm have several problems, 
specially in Steps 3 and 4.
They involve exact
computation of critical points, genericity condition
tests, adjacency tests
and assume exact input equations 
relying on the analysis of the curve at the critical 
points and singularities.
From an algebraic point of view,
they compute sub-resultants
and/or its real roots.
We prove that Steps 3 and 4 can be done
trivially by knowing which point is critical.
In fact, we can make Steps 1,2 and 3 using a subdivision
technique altogether, see Remark \ref{in-ar}.

\subsection*{Existing work.}
Let us mention some of the existing work
and their implementations.
The algorithm mentioned above came up from the CAD method, \citep{MR764184,MR764185}. 
The first implementation that treats singular curves
appears in \citep{MR1116151}.
Steps 3 and 4 of the
algorithm was solved in polynomial time using Thom's codes.
There are dozens of articles improving this
implementation. For example, in \citep{MR1430842,MR3239918}
the authors use bounding boxes to count half-branches
(computing Sturm sequences or counting sign variations).
Also in \citep{MR1008535}
the authors use rational Puiseux expansions 
and in \citep{MR1422725} the authors use
Sturm-Habicht sequences to produce Thom's codes.
The article \citep{MR1106415} compiles 
different implementations and complexities.

\subsection*{Main result.}
In this article we prove that Steps 3 and 4 can be done
by knowing which point is regular and which
point is critical (i.e. singular or has a vertical tangent line),

\begin{nteo}[Main result]
If Steps 0,1 and 2 are done, then Steps 3 and 4 
can be made by $\mathcal{O}(d^3)$ extra additions of integer, where $d$ is the degree of the curve.
\end{nteo}
\begin{dem}
See Theorem \ref{algo} and Theorem \ref{graph}.
\qed
\end{dem}

The basic idea of the proof is that the number of half-branches leaving
a fiber is the same as the number of half-branches
entering the next, in symbols
\[
\sum_{k=1}^dR_{i,k}=\sum_{k=1}^dL_{i+1,k}.
\]
Given that the number of half-branches
on the left and on the right at a regular point is $1$,
with this global invariant, we can determine the number
of half-branches on the left (and on the right) at a 
critical point. We only need to ensure,
that there is at most one critical point per fiber (Step 0).

We first implemented our method using floating point arithmetic
in SageMath \citep{sage}. We noticed that it was very fast, but 
the arithmetic produced several errors. Then,
we implemented our method using polynomials in Bernstein form over the rationals.
It was a slower method but produced a certificated plane curve.
We think that both implementations are useful. The code in SageMath
using floating point arithmetic and some drawings can be found at \url{http://pc-quallbrunn.dm.uba.ar/~cesar/curve2d}.

\subsection*{Summary.}

In the first section we state our main assumption and give
the definition of a \emph{data} for a plane curve. It is basically the information of Steps 1,2 and 3. 
In Sections \ref{data-gr} and \ref{part-data}
we prove that it is possible to construct a graph isotopic to a plane curve $\mathcal{C}$ using Steps $1$ and $2$ without affecting
the total complexity of the algorithm. In Section \ref{refin}
we introduce the notion of \emph{refinement}. This notion
is useful to add more vertexes to the graph in order to produce
a more accurate visualization of the curve.
Finally, in Section \ref{complexity} we count
the operations needed to perform our implementation.

\section{Preliminaries.}
First, let us state our main assumption.

\begin{ass}\label{hypo}
Let $g\in\mathbb{R}[x,y]$ be such that 
the curve $\mathcal{C}=\{g=0\}$ is generically smooth.
Let $\pi_1:\mathbb{R}^2\rightarrow\mathbb{R}$ be the projection to the first coordinate, $\pi_1(x,y)=x$ and 
let $\pi=\pi_1|_\mathcal{C}$. Assume,
\begin{itemize}
\item The map $\pi$ has finite fibers.
\item Every critical point of $\pi$ is in a different fiber. 
\end{itemize}

In concrete examples, we may 
relax the previous assumptions
by asking that every singularity of $\mathcal{C}$ is in a different
fiber of $\pi$. In any case, these assumptions are obtained by a generic change of coordinates.  In
\cite[\S 3.2]{MR1422725}, the authors gave an algorithm to find a change
of variables to produce a curve in general position,
but, as it is said in \cite[\S 3.1]{MR1940259}, 
it is quite expensive even for moderate degrees
in terms of computing times. 
\qed
\end{ass}

\begin{defs}\label{half-branches}
Let $\mathcal{C}$ and $\pi$ be as in Assumption \ref{hypo}.
Let $\{x_i\}_{i=1}^n\subseteq\mathbb{R}$ be an ordered finite subset of 
$\mathbb{R}$. We say that $\{x_i\}_{i=1}^n$ is a \emph{good partition}
(for $\mathcal{C}$ or with respect to $\pi$) if $n\geq 2$ and for every critical value $x$ of $\pi$, 
there exists $i$ such that $x=x_i$ 
and $x_{i-1},x_{i+1}$ are regular values of $\pi$.
 
Clearly, if $\deg(\mathcal{C})=d$, then the number $k$ of critical values
is $\leq d(d-1)$, hence, $n\geq 2k+1$.

\

Let $\{x_i\}_{i=1}^n$ be a good partition for $\mathcal{C}=\{g=0\}$.
For each $1\leq i\leq n$, consider the (possible empty) fiber,
\[
\pi^{-1}(x_i)=\{y_{i1},\ldots,y_{ir_i}\},\quad y_{i1}<\ldots<y_{ir_i},\, 0\leq r_i\leq d.
\]
Let us use notations from \cite[\S 2]{MR1106415}.
Denote by $\rho(x)$ the number of roots of $g(x,-)$, $\rho(x)=\#\pi^{-1}(x)$. Note that $\rho(x)$ remains constant as long as $x$ lies in
a fixed interval $(x_i,x_{i+1})$.
Moreover, 
the corresponding roots are given by continuos
functions 
\[
\varphi_{i,1}<\ldots<\varphi_{i,\rho_i}:(x_i,x_{i+1})\rightarrow\mathbb{R},\
\qquad 1\leq i\leq n-1,\,\rho_i=\max(r_i,r_{i+1}).
\]
The graphs of such functions $\varphi_{i,j}$ are 
the (real) branches of $\mathcal{C}$ over $(x_i,x_{i+1})$.
Note that the continuous extension of 
these functions to $[x_i,x_{i+1}]$ (denoted again as $\varphi_{i,j}$)
satisfies 
$\varphi_{i,j}(x_i)=y_{i,k}$ and 
$\varphi_{i,j}(x_{i+1})=y_{i+1,l}$ for some $k$ and $l$,
$1\leq k\leq r_i$, $1\leq l\leq r_{i+1}$.

The cardinal of the set $\{k\colon\varphi_{ik}(x_i)=y_{ij}\}$
is called \emph{the number of half-branches to the right} 
at $(x_i,y_{ij})$ and $\#\{k\colon\varphi_{i-1,k}(x_i)=y_{ij}\}$
is called \emph{the number of half-branches to the left}
at $(x_i,y_{ij})$.

\

We say that $(\{x_i\},\{y_{ij}\},L,R)$ is a \emph{data}
for $\mathcal{C}$ if $L,R\in\mathbb{N}_0^{n\times d}$,
$L_{ij}$ is the number of half-branches to 
the left and  $R_{ij}$ is the number of half-branches 
to the right at $(x_i,y_{ij})$,
$1\leq i\leq n$, $1\leq j\leq r_i$. If $j>r_i$, we set $L_{ij}=R_{ij}=0$.

The matrices $L$ and $R$ are called
\emph{matrices associated} to $\mathcal{C}$. They are well 
defined because $\pi$ has finite fibers. Also, $L$ and $R$
clearly satisfies
\[
\sum_{k=1}^d R_{i,k}=\sum_{k=1}^d L_{i+1,k},\qquad 1\leq i\leq n-1,
\]
\qed
\end{defs}

\begin{ejemplo}
Consider the nodal plane curve $\mathcal{C}$ given by $g=y^2-x^3-x^2$.
The roots of the $y$-discriminant of $g$ are $\{-1,0\}$.
Hence, we can take a good partition for $\mathcal{C}$
as $\{-2,-1,-0.5,0,1\}$. The fibers are,
\begin{align*}
\pi^{-1}(-2)&=\{g(-2,y)=0\}=\emptyset,\\
\pi^{-1}(-1)&=\{g(-1,y)=0\}=\{0\},\\
\pi^{-1}(-\frac{1}{2})&=\{g(-\frac{1}{2},y)=0\}=
\left\{-\frac{1}{\sqrt{8}},\frac{1}{\sqrt{8}}\right\},\\
\pi^{-1}(0)&=\{g(0,y)=0\}=\{0\},\\
\pi^{-1}(1)&=\{g(1,y)=0\}=\{-\sqrt{2},\sqrt{2}\}.
\end{align*}
Finally, the transposes of $L$ and $R$ are,
\[
L^t=\begin{pmatrix}
0&0&1&2&1\\
0&0&1&0&1
\end{pmatrix},\qquad
R^t=\begin{pmatrix}
0&2&1&2&1\\
0&0&1&0&1
\end{pmatrix}.
\]
\qed
\end{ejemplo}

Finally, let us recall a definition from \citep{Knuth:1976:BOB:1008328.1008329} that we are going
to use in the next sections.
\begin{defs}
Let $n,m:\mathbb{N}\rightarrow\mathbb{N}$ be two functions.
We say that $n$ is \emph{bounded above} by $m$, 
denoted $n=\mathcal{O}(m)$,
if there exists $k\in\mathbb{N}$ such that $n\leq km$.
Analogously, 
We say that $n$ is \emph{bounded below} by $m$, 
denoted $n=\Omega(m)$, if $m$ is bounded above by $n$.

For example, if $\{x_i\}_{i=1}^n$ is a good partition, 
then $n=\Omega(d^2)$.
\qed
\end{defs}

\section{Construction of a graph from a data.}\label{data-gr}

Let us characterize the matrices $L$ and $R$.

\begin{lema}\label{mat-graph}
Given matrices $L,R\in\mathbb{N}_0^{n\times d}$ such that 
\[
\sum_{k=1}^d R_{i,k}=\sum_{k=1}^d L_{i+1,k},\qquad 1\leq i\leq n-1,
\]
there exists a graph with vertexes $\{(i,j)\colon 1\leq i\leq n, 1\leq j\leq d\}$ whose geometric realization has associated matrices $L$ and $R$. Meaning that the point $(i,j)$ is connected by $L_{ij}$ segments to points in $\{(i-1,k)\}_{k=1}^d$ and by $R_{ij}$ segments to points in 
$\{(i+1,k)\}_{k=1}^d$. The intersection of two segments is
included in $\{(i,j)\colon 1\leq i\leq n, 1\leq j\leq d\}\subseteq\mathbb{R}^2$.
\end{lema}
\begin{dem}
We are going to construct the graph iteratively.
If $R_{11}>0$, connect the point $(1,1)$ with the point $(2,j)$,
where $j$ is the first index such that $L_{2j}>0$, that is, $L_{21}=\ldots=L_{2,j-1}=0$.
Then, redefine $R_{11}$ and $L_{2j}$ as $R_{11}-1$ and $L_{2j}-1$.
Repeat this process until $R_{11}=0$.

Inductively, we may assume that $R_{11}=\ldots=R_{1d}=R_{21}=\ldots=R_{i-1,d}=R_{i1}=\ldots=R_{i,j-1}=0$ and $R_{ij}>0$. Again, connect $(i,j)$ with $(i+1,k)$, where $k$ is the first index such that
$L_{i+1,k}>0$. Redefine $R_{ij}$ and $L_{i+1,k}$ as $R_{ij}-1$ and $L_{i+1,k}-1$.
Repeat this process until $R_{ij}=0$.

Finally, we need to check that this process is always possible, and
for this we need the hypothesis. Note that after each redefinition
of $R$ and $L$, the hypothesis is still valid, hence if there exists
some $R_{ij}>0$, there must exists some $L_{i+1,k}>0$.

At the end we obtain a graph with vertexes $(i,j)$ and edges connecting them. Clearly the matrices associated to the geometric realization
of this graph are $L$ and $R$.
\qed
\end{dem}

\begin{prop}\label{graph.constr}
Let $\mathcal{C}=\{g=0\}$ and $\pi$ be as in Assumption \ref{hypo}.
Let $(\{x_i\}_{i=1}^n,\{y_{ij}\},L,R)$ 
be a data for $\mathcal{C}$ and let $d=\deg(g)$.
Then it is possible to construct a graph $G$
from the data that requires $\mathcal{O}(nd)$ additions of integers.
\end{prop}
\begin{dem}
Consider a graph $G$ with vertexes $(x_i,y_{ij})$.
Let $\widehat{G}$ be the graph from Lemma \ref{mat-graph}.
Add an edge to $G$ connecting $(x_i,y_{i,j})$ with $(x_{i+1},y_{i+1,k})$
if there exists an edge in $\widehat{G}$ joining
$(i,j)$ with $(i+1,k)$.
Given that the cardinal of an intermediate fiber of $\pi$ is 
at most $d$, there exist at most $d$ connections between two
fibers, then the total number of connections to construct $G$
is $\mathcal{O}(nd)$. Given that a connection
requires two subtractions, see Lemma \ref{mat-graph},
we obtain that the construction of the graph requires $\mathcal{O}(nd)$ additions of integers.
\qed
\end{dem}

\begin{teo}\label{graph}
Let $\mathcal{C}=\{g=0\}$ and $\pi$ be as in Assumption \ref{hypo}.
Let $(\{x_i\}_{i=1}^n,\{y_{ij}\},L,R)$ 
be a data for $\mathcal{C}$ and let $d=\deg(g)$.
Then the geometric realization 
of the graph $G$ from Proposition \ref{graph.constr}
is isotopic to $\mathcal{C}$
in the region $[x_1,x_n]\times\mathbb{R}$.
\end{teo}
\begin{dem}
Let us use ideas from 
\cite[\S 6.4]{MR2532173} 
to prove that the geometric realization of $G$ 
is isotopic to $\mathcal{C}$. This is a well known result
proved using different methods, see \citep{MR1116151}.

Let us denote by $\varphi_{i,1}<\ldots<\varphi_{i,\rho_i}$
to the (real) branches of $\mathcal{C}$ over $[x_i,x_{i+1}]$, 
see Definition \ref{half-branches}. Also, let 
$s_{i}^{kl}:[x_i,x_{i+1}]\rightarrow\mathbb{R}$
be the linear function whose graph is the geometric
realization of the edge in $G$ with vertexes 
$(x_i,y_{i,k})$ and $(x_{i+1},y_{i+1,l})$,
\[
s_{i}^{kl}(x)=\frac{x_{i+1}-x}{x_{i+1}-x_i}y_{i,k}+
\frac{x-x_{i}}{x_{i+1}-x_{i}}y_{i+1,l}.
\]
Sort the functions $\{s_i^{kl}\}$ 
according to their values in 
$(x_i,x_{i+1})$, $s_{i1}<\ldots<s_{i\rho_i}$.
It is easy to prove (by induction as in Lemma \ref{mat-graph}) that
$\varphi_{ij}(x_i)=s_{ij}(x_i)$ and 
$\varphi_{ij}(x_{i+1})=s_{ij}(x_{i+1})$ for all $i,j$.

\

Let $F:[x_1,x_n]\times\mathbb{R}\times[0,1]\rightarrow [x_1,x_n]\times\mathbb{R}$
be the isotopy defined as $F(x,y,t)=(x,F'(x,y,t))$,
\[
F'(x,y,t)=
(1-t)\cdot y+t\cdot \sum_{i=1}^{n-1} 
f_i(x,y)\cdot\mathbf{1}_{(x_i,x_{i+1})}(x)+
t\cdot y\cdot\mathbf{1}_{\{x_1,\ldots,x_{n}\}}(x)
\]
where $f_i:(x_i,x_{i+1})\times\mathbb{R}\rightarrow\mathbb{R}$
is 
\[
f_i(x,y)=\sum_{j=1}^{\rho_i-1}
\left(\frac{s_{i,j+1}-y}{s_{i,j+1}-s_{i,j}}\varphi_{i,j}+
\frac{y-s_{i,j}}{s_{i,j+1}-s_{i,j}}\varphi_{i,j+1}\right)\cdot 
\mathbf{1}_{(s_{i,j},s_{i,j+1}]}(y)\,+
\]
\[
(y+\varphi_{i,\rho_i}-s_{i,\rho_i}) \cdot\mathbf{1}_{(s_{i,\rho_i},+\infty)}(y)+\,
(y+\varphi_{i,1}-s_{i,1}) \cdot\mathbf{1}_{(-\infty,s_{i,1}]}(y).
\]
Clearly, $F(-,0)$ is the identity, $F(-,1)$ maps the geometric realization
of $G$ to $\mathcal{C}$ and $F(-,t)$ is an embedding. 
The continuity of $F$ follows from the fact that the continuous extension
of $f_i$ to $[x_i,x_{i+1}]\times\mathbb{R}$
satisfies $f_i(x_i,y)=f_i(x_{i+1},y)=y$.
\qed
\end{dem}

\section{Construction of a data from a good partition.}\label{part-data}

Let us give an algorithm to construct the matrices $L$ and $R$
from infinitesimal information attached to regular points of 
$\mathcal{C}$.
This has the advantage of bypassing some computationally expensive algorithms, see \citep{MR1051216,MR1106415,MR3239918}.

\begin{teo}\label{algo}
Let $\mathcal{C}=\{g=0\}$ and $\pi$ as in Assumption \ref{hypo}.
Let $(\{x_i\},\{y_{ij}\})$ be part of a data for $\mathcal{C}$ and
let $g_y$ be the derivative of $g$ with respect to $y$. Let $d=\deg(g)$
and for each $(x_i,y_{ij})$, define $\delta_{ij}\in\{0,1\}$ as 
$\delta_{ij}=1$ iff $g_y(x_i,y_{ij})\neq 0$.
Then, making $\mathcal{O}(d^3)$ extra additions
of integers, it is possible to complete
$(\{x_i\},\{y_{ij}\})$ to a data.
\end{teo}
\begin{dem}
Initialize the matrices $L$ and $R$ according to the following
criteria
\[
L_{ij}=R_{ij}=
\begin{cases}
1&\text{if }g_y(x_i,y_{ij})\neq 0\\
0&\text{if not}
\end{cases}
\]

Now, let us add the information about the critical points of $\pi$.
Recall that 
if $x_i$ is a critical value, then $x_{i-1}$ and $x_{i+1}$ are
regular. Also, that there is at most one critical point per fiber.
For each critical point $(x_i,y_{ij})$, define 
\[
r:=\sum_{k=1}^d L_{i+1,k}-\sum_{k=1}^d R_{i,k},\qquad
l:=\sum_{k=1}^d R_{i-1,k}-\sum_{k=1}^d L_{i,k}
\]
and set $R_{ij}=r$ and $L_{ij}=l$. 
By construction, $(\{x_i\},\{y_{ij}\},L,R)$ is a data
for $\mathcal{C}$.

Then, the number of additions of integers to fill the matrices $L,R$
is $4d$ times the number of critical points.
By B\'ezout, 
\cite[Cor. 18.5]{MR1182558}, 
the set of critical points $\{g=g_y=0\}$ has $\mathcal{O}(d^2)$ points.
Hence, we needed $\mathcal{O}(d^3)$ additions of integers
to construct the matrices $L$ and $R$ from the values $\delta_{ij}$.
\qed
\end{dem}

\begin{obs}
If we only have one singularity per fiber,
for every 
smooth critical point $(x_i,y_{ij})$ we need to 
initialize $(L_{ij},R_{ij})$ as $(2,0)$, $(0,2)$ or $(1,1)$.
According to the next Proposition \ref{crit-poi} this information
can be obtained from the first non-zero value 
of $g_y(x_i,y_{ij}),\ldots,g_{y^d}(x_i,y_{ij})$.
This process can be computationally expensive, see
\cite[Prop. 2.3]{MR1422725},
but it is useful
in practice to avoid, in some cases, the need to rotate the curve.
For example, the curve 
$3 x^{4} + 5 x^{2} y^{2} + 2 y^{4} - 4 y^{2}$, 
an ``$8$'' figure, can be plotted without rotating it to an
``$\infty$'' figure.
\qed
\end{obs}

\begin{prop}\label{crit-poi}
Let $\mathcal{C}=\{g=0\}$ be as in Assumption \ref{hypo}.
and let $(x_0,y_0)\in \mathcal{C}$ be a point such 
that $g_x(x_0,y_0)\neq 0$,
$g_y(x_0,y_0)=g_{yy}(x_0,y_0)=\ldots=g_{y^{r-1}}(x_0,y_0)=0$ and $g_{y^r}(x_0,y_0)\neq 0$ for some $r>0$.

If $r$ is odd, $(x_0,y_0)$ has one half-branch
on the left and one half-branch on the right.
If $r$ is even, we have two possibilities.
If $g_{y^r}(x_0,y_0)g_x(x_0,y_0)>0$, then $(x_0,y_0)$ has two half-branches on the left and if $g_{y^r}(x_0,y_0)g_x(x_0,y_0)<0$, then
$(x_0,y_0)$ has two half-branches on the right.
\end{prop}

\begin{dem}
Consider a parametric curve $\gamma(t)$ given 
by $\gamma(t)=(x(t),t)$, where $x'(0)=0$.
In this special case, we can infer the positions of the half-branches at $\gamma(0)$ using the second derivative test, 
\cite[\S 4.17, Th.4.9]{MR0214705}. 

Note that the curve $\widetilde{\gamma}(t)=(t,x(t))$ parametrizes
the graph of the function $x(t)$. If $x'(0)=0$ and $x''(0)>0$, 
then $x(t)$ has a local minimum at $0$.
This implies that $\gamma$ has two half-branches on the right at $t=0$. Analogously, if $x''(0)<0$, $\gamma$ has two half-branches to
the left at $t=0$. Now, if $x''(0)=0$, we can test $x'''(0)$ and so on.

Assume that $x^{(r-1)}(0)=0$ and $x^{(r)}(0)\neq 0$.
If $r$ is odd, then $\gamma(0)$ has one half-branch on the left
and one half-branch on the right (a saddle-point). If $r$ is even and $x^{(r)}>0$ (resp. $x^{(r)}<0$), then $\gamma(0)$ has two half-branches on the left (resp. on the right).

Now, given that $g_x(x_0,y_0)\neq 0$ and $g_y(x_0,y_0)=0$
we can apply the Implicit Function Theorem, 
\cite[Th.2-12, p.41]{MR0209411}.
Then, there exists a parametric curve $\gamma(t)=(x(t),y_0+t)$,
$|t|<\varepsilon$ such that $\gamma(t)\in \mathcal{C}$,
$\gamma(0)=(x_0,y_0)$ and $\gamma'(0)\neq 0$.
Even more so, $x'(t)=-g_y(\gamma(t))/g_x(\gamma(t))$.
Let us use the following notation:
\[x=x(t),\quad g_y=g_y(\gamma(t)),\quad g_x=g_x(\gamma(t)),\quad
g_{xy}=g_{xy}(\gamma(t)),\quad\ldots.
\]
Let us prove by induction that there exist
rational functions $p_1,\ldots,p_r$ in $t$ and $x(t)$ with denominators of the form $(g_x(\gamma(t)))^k$, $k\in\mathbb{N}_0$,  such that 
\[
x^{(r)}=p_1g_y+\ldots+p_{r}g_{y^{r}},
\] 
where $p_r=-1/g_x$.
The case $r=1$ follows from the Implicit Function Theorem. Assume the case $r$ and let us prove $r+1$.
\[
g_{y^{i}}(\gamma(t))'=g_{y^{i}}'=g_{y^ix}x'+g_{y^{i+1}}
=\frac{-g_{y^ix}}{g_x}g_y+g_{y^{i+1}}
,\quad\forall i\geq 1\Longrightarrow
\]
\[
x^{(r+1)}=(x^{(r)})'=\sum_{i=1}^{r-1}p_i'g_{y^i}+p_ig_{y^i}'
-\frac{g_{y^r}'g_x-g_{y^r}g_{x}'}{g_x^2}=
\sum_{i=1}^{r}\widetilde{q}_ig_{y^i}-\frac{g_{y^{r+1}}}{g_x}.
\]
Now that we have proved the result, 
note that $x'(0)=\ldots=x^{(r)}(0)=0$ if and only if
$g_y(x_0,y_0)=\ldots=g_{y^r}(x_0,y_0)=0$.
Also, if this is the case, $x^{(r+1)}(0)=-g_{y^{r+1}}(x_0,y_0)/g_x(x_0,y_0)$.
\qed
\end{dem}

\section{Refinement.}\label{refin}

\begin{defs}
Let $\mathcal{C}$ and $\pi$ as in Assumption \ref{hypo}.
Let $(\{x_i\}, \{y_{ij}\}, L,R)$ be 
a data for $\mathcal{C}$
and let 
$x\not\in\{x_1,\ldots,x_n\}$. 
We want to add $x$ and
its fiber, $y_1<\ldots<y_s$, to this data.

Given that $x$ is necessarily a regular value of $\pi$, 
we need to insert a column filled with $s$ ones and $d-s$ zeros
to $L$ and to $R$, $(1_1,\ldots,1_s,0_1,\ldots,0_{d-s})$.

If $x<x_1$, we insert the column on the left of $L$ and $R$, if $x>x_n$, we insert the column on the right of $L$ and $R$ and if $x_i<x<x_{i+1}$, 
we insert the column between columns $i$ and $i+1$ of $L$ and $R$.

We say that the addition of a finite set to the data
$(\{x_i\}, \{y_{ij}\}, L,R)$ is a \emph{refinement} of the data.
\qed
\end{defs}

\begin{obs}
The process of refinement is important in the algorithms 
for spacial curves and surfaces to produce a 
good partition common to all the plane curves
involved. For example, in the spacial curve case, we take a data
for the projection of the curve to the $xy$-plane and a data
for the projection to the $xz$-plane. Then, we \emph{refine} both
partitions to produce a common one.

Another application of this process is the possibility of
adding more points to a plane curve without affecting the general
complexity of the algorithm, we only need to find
real roots of univariate polynomials, see Lemma \ref{fp}.
\qed
\end{obs}

\section{Complexity.}\label{complexity}

In the next lemma, we only pay attention to
the computations that need floating point arithmetic, 
that is, finding real roots and making additions
and multiplications.
Derivatives, discriminants and resultants 
can be done symbolically. For the complexity
of real root isolation, see \citep{MR3424039,MR2474343}.

\begin{lema}[\textbf{Floating point arithmetic}]\label{fp}
Let $\mathcal{C}=\{g=0\}\subseteq\mathbb{R}^2$ be a 
degree $d$ plane curve as in Assumption \ref{hypo}.
Let $\{x_i\}_{i=1}^n$ 
be a good partition for $\mathcal{C}$. Recall $n=\Omega(d^2)$.

Then, in order to complete $\{x_i\}_{i=1}^n$ to a data it is sufficient
to find the real roots of $n$ degree $d$ univariate polynomials
and to make $\mathcal{O}(nd^2)$ multiplications
and additions.

Even more, the refinement of a data by some $x$, 
requires finding the real roots of one degree $d$ univariate 
polynomials and $\mathcal{O}(d^2)$ multiplications
and additions.
\end{lema}
\begin{dem}
In order to compute a fiber $\pi^{-1}(x)$, 
we need to compute $\{y_j\}_{j=1}^s$, the real roots
of a degree $d$ univariate polynomial $g(x,-)$.
Also, we need to evaluate $g_y$ at every point $(x,y_j)$, $1\leq j\leq s\leq d$.
As we proved in Theorem \ref{algo}, 
with these information we can compute a data (and a refinement)
for $\mathcal{C}$.

Using Horner's method, \cite[\S 4.6.4]{MR3077153}, the evaluation of $g_y$ at $(x,y_j)$ requires 
$\mathcal{O}(d^2)$ operations.
\qed
\end{dem}

\begin{obs}
If we use Proposition \ref{crit-poi} to complete a good partition to a data (as in Lemma \ref{fp}), we need to 
evaluate $g_y(x_i,y_{ij}),\ldots,g_{y^d}(x_i,y_{ij})$
at every critical point. This requires $\mathcal{O}(d^5)$ extra operations. Thus, if we start with a good partition
with $n=\Omega(d^3)$ elements or if we refine the data
with $\Omega(d^3)$ elements, then 
the computations using Proposition \ref{crit-poi}
does not alter the complexity of the algorithm.
\qed
\end{obs}

\begin{obs}[Interval arithmetic]\label{in-ar}
The success of our algorithm relies on finding real roots accurately.
Standard numerical methods do not produce a certificated answer.
We propose, as in \citep{MR2381611}
to work with polynomials in Bernstein form (see \citep{MR2499511,MR3079719})
and a bisection method. Also, we assume that $g\in\mathbb{Q}[x,y]$.
First, we find disjoint intervals
of length $\varepsilon$ containing the roots of 
the discriminant $\Delta(g)$, $\varepsilon\in\mathbb{Q}$.
Then, adding more intervals, we define a 
good partition $\{[x_i,x_i+\varepsilon]\}_{i=1}^n$, $x_i\in\mathbb{Q}$.
Second, using a Bernstein basis in two variables, 
find disjoint intervals of length $\varepsilon$, $[y_{ij},y_{ij}+\varepsilon]$, such that $\pi^{-1}([x_i,x_i+\varepsilon])$ is included in $\bigcup_{j} [x_i,x_i+\varepsilon]\times [y_{ij},y_{ij}+\varepsilon]$
, $y_{ij}\in\mathbb{Q}$.
Finally, complete $(\{[x_i,x_i+\varepsilon]\},\{[y_{ij},y_{ij}+\varepsilon]\})$ to a data. In order to do so, we need the sign of $g_y$ at each rectangle $[x_i,x_i+\varepsilon]\times [y_{ij},y_{ij}+\varepsilon]$, and this can be done effectively writing $g_y$ in  Bernstein form.

Note that it is possible to avoid the computation of $\Delta(g)$ finding 
rectangles $[x,x+\varepsilon]\times[y,y+\varepsilon]$ that contain the solutions of $g=g_y=0$.
\qed
\end{obs}

%

\bibliographystyle{elsarticle-harv}
\bibliography{citas}
\end{document}